\documentclass[a4paper, twoside, 12pt]{article}

\usepackage[margin=3cm, headheight=14.5pt]{geometry} 
\usepackage{amsmath, amssymb, amsthm} 
\usepackage{graphicx} 
\usepackage{hyperref} 
\usepackage{listings} 
\usepackage{fancyhdr} 
\usepackage{titling} 
\usepackage{titlesec} 
\usepackage{tcolorbox}
\usepackage{tikz}
\usepackage{mathtools}
\usepackage{arydshln}
\usepackage{bbold}
\usepackage{titling}
\usepackage{lipsum} 
\usepackage{ragged2e}
\usepackage{titlesec}
\usepackage{sectsty}
\usepackage{wasysym}
\usepackage{mathtools}
\usepackage{ MnSymbol}
\usepackage[all,cmtip]{xy}
\usepackage{ytableau}
\usepackage{enumitem}
\usepackage{tocloft}
\usepackage[T1]{fontenc}
\usepackage{imakeidx}
\usepackage{braket}
\usepackage{setspace}

\sectionfont{\large\bfseries}
\subsectionfont{\normalsize\bfseries}
\subsubsectionfont{\small\bfseries}

\makeatletter
\renewcommand{\l@section}[2]{%
  \ifnum \c@tocdepth >\m@ne
    \addpenalty{-\@highpenalty}%
    \vskip 1.0em \@plus\p@
    \setlength\@tempdima{1.5em}%
    \begingroup
      \parindent \z@ \rightskip \@pnumwidth
      \parfillskip -\@pnumwidth
      \leavevmode \bfseries
      \advance\leftskip\@tempdima
      \hskip -\leftskip
      #1\nobreak\hfil \nobreak\hb@xt@\@pnumwidth{\hss #2}\par
    \endgroup
  \fi}
\makeatother

\theoremstyle{plain}
\newtheorem{theorem}{Theorem}[section]
\newtheorem{lemma}[theorem]{Lemma}
\newtheorem{corollary}[theorem]{Corollary}
\theoremstyle{definition}
\newtheorem{definition}[theorem]{Definition}
\newtheorem{example}[theorem]{Example}

\newtheorem{remark}[theorem]{Remark}
\newtheorem{proposition}[theorem]{Proposition}

\numberwithin{equation}{section}

\usepackage{titlesec}
\usepackage{sectsty}
\usepackage{indentfirst,csquotes}
\usepackage{csquotes}
\usepackage[style=numeric,backend=biber]{biblatex}
\titleformat{\section}[hang]{\Large\bfseries\raggedright}{\thesection}{1em}{} 
\titleformat{\subsection}[hang]{\large\bfseries\raggedright}{\thesubsection}{1em}{} 
\subsubsectionfont{\normalsize}
\titleformat{\subsubsection}[runin]{\normalsize\bfseries\raggedright}{\thesubsubsection}{1em}{}

\makeindex

\begin{document}


\thispagestyle{empty}
\begin{center}
{\Large  \MakeUppercase{\bf Localisation in Equivariant Cohomology}\\}

\vspace*{1cm}

\noindent
\begin{tabular}[t]{@{}c}
\textbf{Catherine C. Notman}\\  MPhys Mathematical Physics \\ University of Edinburgh \\ C.Notman@sms.ed.ac.uk
\end{tabular}
\hfill
\begin{tabular}[t]{@{}c}
\textbf{Muaadh A. Sanabani} \\ MPhys Theoretical Physics \\ University of Edinburgh \\ M.Sanabani@sms.ed.ac.uk
\end{tabular}

\vspace*{0.5cm}
{\today}
\begin{abstract}
Equivariant cohomology, a captivating fusion of symmetry and abstract mathematics, illuminates the profound role of group actions in shaping geometric structures. At its core lies the Atiyah-Bott Localization Theorem, a mathematical jewel unveiling the art of localization. This theorem simplifies intricate integrals on symplectic manifolds with Lie group actions, revealing the hidden elegance within complexity. Our paper embarks on a journey to explore the theoretical foundations and practical applications of equivariant cohomology, demonstrating its transformative power in diverse fields, from theoretical physics to geometry. As we delve into the symphonic interplay between geometry and symmetry, readers are invited to witness the beauty of mathematical patterns emerging from abstraction. This mathematical voyage unveils the harmonious marriage of symmetry, topology, and elegance in the captivating realm of equivariant cohomology.
\end{abstract}

\end{center}
\tableofcontents
\thispagestyle{empty}
\newpage



\pagenumbering{arabic}

\raggedright 
\justify

\section{Chapter 1: Introduction}

\subsection{Background}
The captivating interplay between symmetry and mathematical structures has been a driving force in both theoretical mathematics and physics for centuries. Symmetry not only enhances our understanding of complex phenomena but also enriches our ability to unravel the elegant patterns underlying these intricate systems. Within the realm of algebraic topology and geometry, the field of equivariant cohomology emerges as a profound mathematical framework that artfully marries symmetry and mathematical abstraction. \newline

\noindent Equivariant cohomology serves as a powerful extension of classical cohomology theory, ingeniously designed to accommodate the symmetries introduced by group actions on topological spaces. This remarkable branch of mathematics equips us with the tools to dissect, analyze, and quantify the effects of symmetries on geometric structures, paving the way for deeper insights into the fundamental nature of these structures. \newline

\noindent In this pursuit, our exploration takes us on a captivating journey through the world of equivariant cohomology, where symmetries become the guiding threads in the tapestry of mathematical investigation. At the heart of our quest lies the profound Atiyah-Bott Localization Theorem, a mathematical gem that simplifies intricate integrals over symplectic manifolds endowed with Hamiltonian actions of Lie groups. Through this theorem, we uncover the art of localization—a technique that dissects complex integrals into contributions from isolated fixed points, revealing a hidden elegance in the midst of mathematical complexity.\newline

\noindent This paper embarks on a quest to decipher the intricacies of equivariant cohomology, unveiling its theoretical foundations, applications, and profound implications. We will explore the central role of symmetry and group actions in shaping equivariant cohomology and dive into the mechanics of localization techniques that transform daunting integrals into elegant sums over fixed points. Through illustrative examples and practical applications, we will demonstrate the practical utility of these mathematical tools in diverse fields, including theoretical physics and geometry.\newline

\noindent As we journey through the landscape of equivariant cohomology and localization, we invite readers to embrace the beauty of mathematical symmetries and their profound impact on our understanding of the world. The symphony of geometry and symmetry unfolds before us, revealing patterns and structures that transcend mere abstraction. Join us on this mathematical voyage—a voyage that showcases the harmonious marriage of symmetry, topology, and elegance within the captivating realm of equivariant cohomology.\newline

\subsection{Overview of Localisation Theorem}
\begin{figure}[ht]
\begin{minipage}[c]{\textwidth}
\centering
\includegraphics[width=0.50\textwidth]{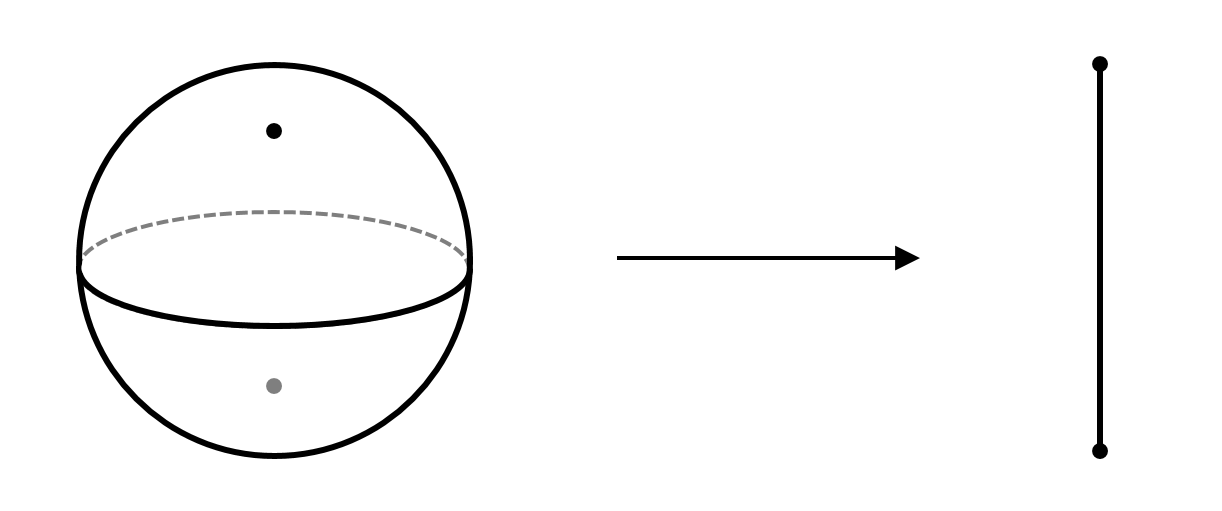}
\caption{The acion of $S^{1}$ on $S^{2}$ and the associated moment map $\mu$ \cite{symplctic}.}
\label{fig_shpere_mot}
\end{minipage}%
\end{figure}

Finding the surface area of a unit sphere $S^{2}$ using the framework of equivariant cohomology provides an overview of the elegance of this method. By exploiting the rotational symmetry of $S^{2}$, we find that the integral of the area form $\omega$ boils down to a summation of two points.
$$
\int_{S^{2}} \sin\phi \, d\theta\wedge d\phi \quad \quad \lbrace 0 \leq \theta \leq 2\pi, \ 0 \leq \phi \leq \pi \rbrace
$$
Let $S^{1} \circlearrowright S^{2}$, the moment map of this action is simply the height function, and the fixed points of the action are the north and south poles of $S^{2}$ (see figure \ref{fig_shpere_mot}).
 \newline
The fundamental vector field $X = - 2\pi \frac{\partial}{\partial \theta} \in \operatorname{Lie}(S^{1})$ provides basis for the algebra of the circle group. By definition of the Hamiltonian action $\iota_{X}\omega = d \mu$,
\begin{align*}
\iota_{X} \omega &= 2\pi \sin \phi \, d\phi\\
&= d(2\pi \cos\phi)
\end{align*}
The moment map $\mu = 2\pi \cos \phi$. Now, let's look at the fixed points of the action, 
$$
F \coloneqq (S^{2})^{S^{1}} = \lbrace \phi=0, \phi = \pi \rbrace 
$$
Applying the localization formula, \emph{theorem} \ref{localisation_form}, we find the surface area of the unit sphere 
$$
\int_{S^{2}} \sin\phi \, d\theta\wedge d\phi = 2\pi - (-2\pi) = 4\pi \quad \quad \square
$$


\newpage

\section{Chapter 2: Mathematical Preliminaries}
\subsection{Topological Manifolds}

\definition {\it{Topological Manifolds}} \cite{NAK}. A topological manifold $M$ is a topological Hausdorff space that has a countable basis for its topology and that is locally homeomorphic to $\mathbb{R}^n$. The number $n$ is called the dimension of $M$. 

\begin{figure}[ht]
\begin{minipage}[c]{\textwidth}
\centering
\includegraphics[width=0.6\textwidth]{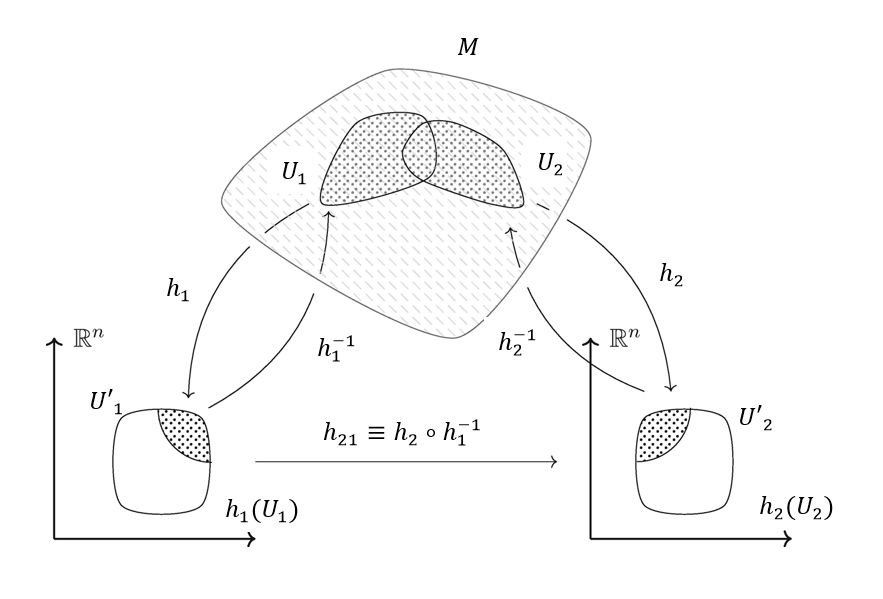}
\caption{Charts and transition map on an $n$-dimensional manifold $M$ \cite{Nima}.}
\label{fig_smooth_mapping}
\end{minipage}%
\end{figure}

\definition {\it{Smooth Structure}} \cite{From_Calc} on a topological manifold,

\begin{enumerate}[label=\roman*]
    \item A chart $(U,h)$ on an $n$-dimensional manifold is a homeomorphism $h: U \rightarrow U'$, where $U$ is an open set in $M$ and $U'$ is an open set in $\mathbb{R}^n$.
    \item A system $\mathcal{A} = \lbrace h_i: U_i \rightarrow U'_i \ | \ i\in J \rbrace$ of charts is called an atlas, provided $\lbrace U_i\ | \ i \in J \rbrace$ covers $M$. 
    \item An atlas is smooth when all the maps
    $$
    h_{ji} = h_j \circ h^{-1}_i : h_i(U_i \cap U_j) \rightarrow h_j(U_i \cap U_j)
    $$
    are smooth. These are called {\it{transition functions}} for the atlas (see figure \ref{fig_smooth_mapping}).
\end{enumerate}

\definition {\it{Smooth Manifold}} \cite{Sch2017}. A smooth manifold is a pair $(M,\mathcal{A})$ consisting of a topological manifold $M$ and a smooth structure $\mathcal{A}$ on $M$.\newline 

\noindent Usually $\mathcal{A}$ is suppressed from the notation, we will write $M$ instead of $(M,\mathcal{A})$.

\example The $n$-dimensional sphere $S^n = \lbrace x \in \mathbb{R}^{n+1} \ | \  |x| = 1 \rbrace$ is an $n$-dimensional smooth manifold. We define an atlas with $2(n+1)$ charts $(U_{\pm i}, h_{\pm i})$ where
$$
U_{+i} = \lbrace  x \in S^n \ | \ x_i > 0 \rbrace, \quad\quad\quad U_{-i} = \lbrace  x \in S^n \ | \ x_i < 0 \rbrace
$$
and $h_{\pm i}: U_{\pm i} \rightarrow D^n$ (where $D^n$ is an open disc in $\mathbb{R}^n$) is the map given by $h_{\pm i}(x) = (x_1, \cdots, \hat{x}_i, \cdots, x_{n+1})$. The circumflex over $x_i$ denotes that $x_i$ is omitted. The inverse map is
$$
h_{\pm i}^{-1} = \left( u_1, \cdots, u_{i-1}, \pm \sqrt{1- |x|^2}, u_i, \cdots, u_n \right)
$$

\begin{figure}[ht]
\begin{minipage}[c]{\textwidth}
\centering
\includegraphics[width=0.7\textwidth]{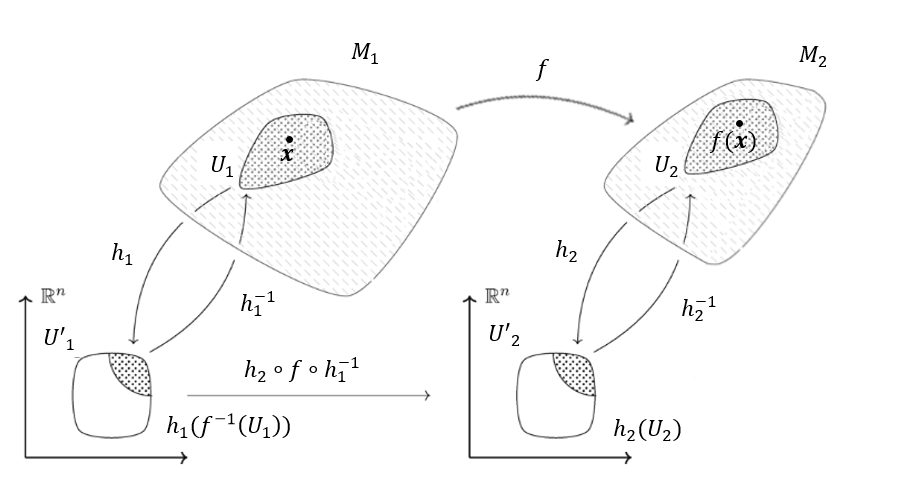}
\caption{Smooth map between manifolds $M_1$ and $M_2$ \cite{Lee12}.}
\label{fig_diffeo} 
\end{minipage}%
\end{figure}

\definition {\it{Smooth Map}} \cite{From_Calc}. Consider smooth manifolds $M_1$ and $M_2$ and a continuous map $f: M_1 \rightarrow M_2$. The map $f$ is called smooth at $x\in M_1$ if there exists charts $h_1: U_1 \rightarrow U'_1$ and $h_2: U_2 \rightarrow U'_2$ on $M_1$ and $M_2$ with $x \in U_1$ and $f(x) \in U_2$ such that,
$$
h_2 \circ f \circ h^{-1}_1 : h_1(f^{-1}(U_2)) \longrightarrow U'_2
$$
is smooth in a neighbourhood of $h_1(x)$. If $f$ is smooth at all points of $M_1$ then $f$ is said to be a smooth map (see figure \ref{fig_diffeo}).  

\definition {\it{Submanifold}} \cite{LOR08}. A subset $N \subset M^n$ of a smooth manifold is said to be a smooth submanifold (of dimension $k$) if the following condition is satisfied: for every $x \in N$ there exists a chart $h: U \longrightarrow U'$ on $M$ such that 
$$
x \in U \quad \quad \text{and} \quad \quad h(U \cap N) = U' \cap \mathbb{R}^k 
$$
where $\mathbb{R}^k \subseteq \mathbb{R}^n$ is the standard subspace. 

\example The $n$-sphere $S^n$ is a smooth submanifold of $\mathbb{R}^{n+1}$. 

\definition {\it{Embedding}} \cite{From_Calc}. An embedding is a smooth map $f: N \longrightarrow M$ such that $f(N) \subset M $ is a smooth submanifold and $f: N \longrightarrow f(N)$ is a diffeomorphism. 

{\theorem {\normalfont{(Whitney embedding \cite{NAK})}}. Let $M^n$ be a smooth manifold of dimension $n$. There exists an embedding of $M^n$ into a Euclidean space $\mathbb{R}^{n+k}$ for some $k \in \mathbb{N}$.} 

{\corollary Every compact topological $n$- dimensional manifold is homeomorphic to a (locally flat) topological submanifold of a Euclidean space $\mathbb{R}^{n+k}$ for some $k \in \mathbb{N}$.}

\definition {\it{Tangent space}} \cite{NAK}. The tangent space to a point on a smooth $n$-manifold may be thought as the set of tangent vectors to all possible curves passing through that point. Formally it is the set of all derivatives and is spanned by $\lbrace \frac{\partial}{\partial x^1}, \frac{\partial}{\partial x^2}, \ldots \frac{\partial}{\partial x^n} \rbrace $ where $x^1, x^2, \ldots,x^2 $ are local coordinates. 
The cotangent space $T_p^*M$ is the dual space to $T_pM$ and is spanned by the one-forms $\lbrace dx^1, dx^2, \ldots, dx^n \rbrace $ which are defined by their action on the basis vectors of $T_pM$: $ dx^i(\frac{\partial}{\partial x^j}) = \delta_j^i \,. $ \newline

\definition \emph{Differential $k-$forms} \cite{MartensGeo}. A $k-$form at $p\in M$ is a map of $k$ vectors in $T_pM$ to $\mathbb{R}$ which is \emph{multilinear} (linear in each argument) and \emph{alternating} (changes sign under exchange of two arguments). The set of differential $k-$forms on $M$ is denoted by $\Omega^{k}(M)$.

\remark The \emph{differential forms} form a graded algebra.
$$
\Omega(M) = \bigoplus_{i=0}^{\operatorname{dim}M} \Omega^i(M)
$$

\definition {\it The wedge product }\cite{LOR08}. The wedge product $\wedge$ is a bilinear alternating map between two differential forms resulting a higher degree form.
\begin{align*}
\wedge: \Omega^{k}(M) \times \Omega^{n}(M) &\longrightarrow \Omega^{k+n}(M) \\
\wedge (\omega, \alpha) &\mapsto \omega \wedge \alpha
\end{align*}
such that $\omega \wedge \alpha = (-1)^{kn} \alpha  \wedge \omega$  $\forall \omega \in \Omega^{k}(M), \ \alpha \in \Omega^{n}(M)$.

\definition \emph{The exterior derivative }\cite{MartensGeo}. The exterior derivative is an $\mathbb{R}-$linear map that satisfies two conditions. 
$$
d: \Omega^k(M) \rightarrow \Omega^{k+1}(M)
$$
\begin{enumerate}[label=(\roman*)]
\item $d(\omega \wedge \alpha ) = d\omega \wedge \alpha + (-1)^{k} \omega \wedge d\alpha  $
\item $d(d\omega) = 0 \ \forall \omega \in \Omega^{k}(M)$, or more compactly $d^2 = 0$ 
\end{enumerate}
\definition \emph{Exact \& Closed forms} \cite{NAK}. A form $\alpha$ is closed if $d\alpha = 0 $ and exact is there exists a $(k-1)-$form $\phi$ satisfying $\alpha = d\phi$. All top forms $\alpha \in \Omega^{dim(M)}$ are closed. 

\definition \emph{de Rham complex} \cite{LOR20}. The properties of the exterior derivative gives rise to the complex 
$$0 \xrightarrow{d} \Omega^0(M) \xrightarrow{d} \Omega^1(M) \xrightarrow{d} \cdots \xrightarrow{d} \Omega^{n-1}(M) \xrightarrow{d} \Omega^n(M) \xrightarrow{d} 0 \,.$$

\remark Any consecutive pair of links in the chain the kernel of second map is contained in the image of the first; the former is the set of closed $k$-forms, denoted $Z^k$, and the later is the set of exact $k$-forms, denoted $B^k$, both of which are subspaces of $\Omega^k(M)$. 

\definition \emph{de Rham cohomology }\cite{NAK}. The de Rham cohomology studies how much closed forms fail to be exact. The $k$-th de Rham cohomology is the space of equivalences classes of closed $k$-forms that differ by and exact $k$-form:
$$ 
H^k(M) = Z^k(M) / B^k(M) = Z^k(M) / \sim \, = \lbrace [x] : x\sim y \iff x-y \in B^k(M) \rbrace \,.
$$
using the de Rham complex, we can construct and alternative definition of de Rham cohomology 
$$
H^{k}(M) \coloneqq \frac{\operatorname{ker}d: \Omega^{k}\rightarrow \Omega^{k+1}}{\operatorname{im}d: \Omega^{k-1}\rightarrow \Omega^{k}}
$$

\remark The cohomology forms a graded algebra over a ring $R$ $$H(M; R) = \bigoplus_{i=0}^{dim\, M} H^i(M;R) \,.$$

\definition \emph{Integration on manifolds} \cite{LOR20}. Integrating a differential form on a manifold is a linear mapping between the top cohomology of the manifold and the cohomology of a point. 
\begin{align*}
\int : H^{\bullet}(M;R) &\rightarrow H^{\bullet}(pt;R)
\end{align*}

To study the relationship between two different topological manifolds one can inquire as to whether they are homotopic. 
\definition \emph{homotopy class }\cite{LOR20}. Let $X$ and $Y$ be topological spaces and let $x_0$ and $y_0$ respectvely be base points on each. Two maps $f_1, f_2 : (X, x_0) \rightarrow (Y, y_0) $ are \emph{homotopic} if there is exists a continuous map $F : X \times [0, 1] \rightarrow Y$ which satisfies for all $t\in [0, 1]$ 
$$
F(x, 0) = f_1(x) \,,\quad F(x, 1) = f_2(x) \,,\quad F(x_0, t) = y_0 \,.
$$
The last relation means that the map is base point preserving. Homotopy defines an equivalence relation of maps, so that two maps are in the same \emph{homotopy class} if they are homotopic. 

\definition \emph{Homotopy Equivalence }\cite{NAK}. Let $X$ and $Y$ be topological spaces, we say that $X$ and $Y$ are homotopy equivalence if there exists two maps $f: X \rightarrow Y$ and $g: Y \rightarrow X$ such that $f \circ g = \mathrm{id}_{Y}$ and $g \circ f = \mathrm{id}_{X}$.

\definition \emph{The fundamental group }\cite{LOR20}. $\pi_(X, x_0)$ of a topological space $X$ is the set homotopy classes of maps $f : (S^1, s_0) \rightarrow (X, x_0) $ with concatenation of maps as the group operation. 

It is possible to generalize this concept for $q\geq 1$ be defining the $q$-th homotopy group $\pi_q(X, x_0)$ as the set of homotopy classes of maps $f: (S^q, s_0) \rightarrow (X, x_0) \,. $ \newline

To understand the properties of vector spaces and their extensions, it's essential to explore the concept of complexification. We might need to prove certain aspects of this section as they are extremely important for later sections.
 
\definition \emph{Complexification }\cite{complexif}. A complexification of a real vector space $V$ is given by the tensor product of the vector space $V$ and the complex plane $\mathbb{C}$ along the real ring $\mathbb{R}$. 
$$V\oplus_{\mathbb{R}}\mathbb{C}$$ 

{\proposition The complexification of a real vector space $V$ is isomorphic to the direct sum of a complex vector space and its conjugate}
$$
V \otimes_{\mathbb{R}} \mathbb{C} \cong E \oplus \overline{E} 
$$
\proof Let's consider the $\mathbb{R}$ linear map $J: V \rightarrow V$ such that $J  \circ J = -\mathrm{id}$. This gives rise to the complex map $J_{\mathbb{C}}: V \otimes_{\mathbb{R}} \mathbb{C} \rightarrow V \otimes_{\mathbb{R}} \mathbb{C}$ between the complexification of the $\mathbb{R}$ vector space $V$ and itself, defined by multiplying by $\sqrt{-1}$. We know that $J_{\mathbb{C}} \circ J_{\mathbb{C}} = -\mathrm{id}$.\newline

\noindent Since $J_{\mathbb{C}}$ is a complex linear map, we can view $V \otimes_{\mathbb{R}} \mathbb{C}$ as a complex vector space. Thus, $J_{\mathbb{C}}$ has eigenvalues $\pm i$. Let's denote the corresponding eigenspaces as $E(i)$ and $E(-i)$, respectively.\newline

\noindent The eigenspace $E(i)$ consists of all elements $v \in V \otimes_{\mathbb{R}} \mathbb{C}$ such that $J_{\mathbb{C}}(v) = i v$. Similarly, $E(-i)$ consists of all elements $v$ such that $J_{\mathbb{C}}(v) = -i v$.\newline

\noindent We can show that $E(i)$ and $E(-i)$ are complex subspaces of $V \otimes_{\mathbb{R}} \mathbb{C}$. Let $v_1, v_2 \in V \otimes_{\mathbb{R}} \mathbb{C}$ and $c \in \mathbb{C}$. Then we have:
\[J_{\mathbb{C}}(cv_1 + v_2) = cJ_{\mathbb{C}}(v_1) + J_{\mathbb{C}}(v_2) = c(iv_1) + iv_2 = i(cv_1 + v_2),\]
which implies that $cv_1 + v_2 \in V \otimes_{\mathbb{R}} \mathbb{C}$. Similarly, we can show that $\overline{E}$ is closed under scalar multiplication and vector addition.\newline

\noindent Now, we can observe that $V \otimes_{\mathbb{R}} \mathbb{C} \cong E(i) \oplus E(-i),$ i.e., any element $v \in V \otimes_{\mathbb{R}} \mathbb{C}$ can be uniquely decomposed as $v = v_i + v_{-i}$, where $v_i \in E(i)$ and $v_{-i} \in E(-i)$. Therefore we may define an isomorphism $V \otimes_{\mathbb{R}} \mathbb{C} \cong E(i) \oplus E(-i)$ by  $v = v_i + v_{-i} \mapsto (v_i, \, v_{-i}).$ \newline

\definition \emph{Complex structure} \cite{NAK}. Let $M$ be a manifold with real dimension $2n$ spanned by coordinated $\lbrace x_1, ..., x_n, y_1, ... , y_n \rbrace$. As a real vector space $TM$ is spanned by $\lbrace \frac{\partial}{\partial x_1}, \ldots,\frac{\partial}{\partial x_n}, \frac{\partial}{\partial y_1}, \ldots, \frac{\partial}{\partial y_n} \rbrace $. A complex structure is a map $J$ such that $$J\frac{\partial}{\partial x_i} = \frac{\partial}{\partial y_i}, \quad J\frac{\partial}{\partial y_i} = -\frac{\partial}{\partial x_i}  $$ 
where $J^2 = - \mathrm{id}$. Such a structure is compatible with a metric $g$ on $M$ 
$$\forall X, Y \in TM \ g(JX, JY) = g(X, Y)$$

\remark The complexification of the cotangent space at a point $T_p^*M \otimes_{\mathbb{R}} \mathbb{C}$ which has basis $\lbrace dz^j|_p = (dx^j|_p + idy^j|_p), \ d\overline{z}^j|_p = (dx^j|_p - idy^j|_p)  : j = 1, ..., n \rbrace$. Thus $T_pM \otimes_{\mathbb{R}} \mathbb{C}$ has basis $\lbrace \frac{\partial}{\partial z^i}=\frac{1}{2}(\frac{\partial}{\partial x^i}-\frac{\partial}{\partial y^i}), \frac{\partial}{\partial \overline{z}^i}=\frac{1}{2}(\frac{\partial}{\partial x^i}+\frac{\partial}{\partial y^i}) : i = 1, ... , n \rbrace.$ 

\definition {\it Kahler Manifold }\cite{carosso2018geometric}. Let $M$ be a complex manifold with a complex structure $J$ compatible with the Riemannian metric $g$. The \emph{Kahler form} is a second rank tensor $\omega$ defined by $\omega(X, Y) = g(JX, Y)$. If $\omega$ is closed then locally there exists a real function $\mathcal{K}$ called the Kahler potential, such that $\omega= i\partial\overline{\partial}\mathcal{K}$.  Then $(M, J, g, \mathcal{K})$ is a \emph{Kahler manifold}.


\subsection{Basics of Lie Groups and Lie Algebras}
Lie groups are of fundamental importance to the study of mathematical structures and physical systems, as they serve the purpose of encoding symmetries. Symmetries permeate all fields within the study of physics and are essential for the development and construction of theories. They provide essential information about the object of interest and the properties of systems. Within the Lagrangian formalism symmetries of the Lagrangian give rise to conserved quantities, thus revealing vital information about the system.

\definition{\it Lie Groups }\cite{Sch2017}. A Lie group is a smooth manifold $G$ with action $\bullet$ for which the multiplication and inversion maps, $\mu$ and $i$ respectively, are both smooth.
\begin{center}
    \begin{tabular}{cc}
$\mu \, : \, G\times G  \rightarrow G$    &  $i\, :\, G \rightarrow G $  \\
$(\rho, \rho') \mapsto \rho \bullet \rho'$   &  $\rho \mapsto \rho^{-1}$ 
\end{tabular}
\end{center}

\definition{\it Lie Algebras }\cite{Sch2017}. A Lie algebra ($\mathbf{g}$, $\left[\cdot ,\cdot\right]$) is a vector space over the field $K$ with a symmetric bilinear form $[\cdot,\cdot]$, which satisfies the Jacobi identity,

 $$    \forall X, Y, Z \in \mathbf{g} \quad [X, [Y, Z]] + [Z, [X, Y]] + [Y, [Z, X]] = 0 \,. $$

\remark The Lie algebra $\mathbf{g}$ of a Lie Group $G$ is isomorphic as a vector space to the tangent space of $G$ at the identity: $T_eG \cong \mathbf{g}$. \newline

Within this work the Lie groups of interest will be subgroups of the general linear group $GL(n, K)$ for $K = \mathbb{R} \text{  or  } \mathbb{C}$, most often $SU(n)$ or the torus group of dimension $n$, $T^n = U(1)^n$, in addition to the symplectic group $Sp(2n)$.

\definition{\it Right action }\cite{NAK}. The right action of a Lie group on a smooth manifold $M$ is defined by a smooth map $G\times M \rightarrow M$ taking $(\rho, x) \mapsto \rho\bullet x$ such that the action is associative, $\forall \rho, \rho' \in G \,,\, (\rho\bullet \rho') (x) = \rho\bullet \rho'(x) $, and all points of $M$ are invariant under the identity element $e\in G$ such that $e(x) = x.$

\definition{\it The orbit space }\cite{Lie-rep-notes}. The orbit space $M / G$ is a quotient space given by the equivalence relation $x \sim y \implies y \in Orb_G(x).$ The equivalence classes are the orbits of $M$ under the action of $G.$

\remark Let $M$ be a smooth manifold and $x\in M $. The orbit of $x$ under the action of a Lie group $G$ is the set of points of the form $\rho\bullet x$ for $\rho$ in  $G$, which can be expressed as $Orb(x)= \lbrace \rho\bullet x | \rho \in G \rbrace = im_G(x)\subseteq M.$ A fixed point is an element $x$ in $M$ such that its orbit contains only itself; that is, all elements of $G$ act trivially on it. The fixed point set is the set of all such points.

\definition{\it The stabilizer }\cite{Lie-rep-notes}. The stabilizer of a point $x$ is the set of element $\rho\in G$ the leave $x$ invariant; $Stab(x) = \lbrace \rho \in G | \forall x \in M \,, \rho\bullet x = x \rbrace .$ The group $G$ acts freely on $M$ is the stabilizer of every point of $M$ contains only the identity element.

{\theorem {\normalfont{(\cite{maximaltori})}} Any connected Lie group $G$ is isomorphic to $(\mathbb{R}/\mathbb{Z})^m \times \mathbb{R}^k$ for non-negative integers $k$ and $m$. In particular, if $G$ is compact, it is isomorphic to $(\mathbb{R}/\mathbb{Z})^m$ for some $m$, this is defined as the \textbf{maximal torus} of $G$
}

\definition \emph{Torus} \cite{Lie-rep-notes}. A compact connected abelian Lie group $T$ is called a torus by the previous theorem, $T\cong (\mathbb{R}/\mathbb{Z})^m $ for a non-negative integer $m$. 

\remark If $T$ is a torus as in the definition, its Lie algebra $\mathbf{t}$ can be identified with $\mathbb{R}^m$, and its exponential map is given by

$$
exp: \mathbf{t}= (t_1,t_2,\ldots,t_m) \in \mathbb{R}^m \mapsto \mathbf{z}= (z_1,z_2,\ldots,z_m)
$$
where $z_j = e^{2\pi i t_j}$. This way we explicitly describe the quotient map $\mathbb{R}^m \rightarrow (\mathbb{R}/\mathbb{Z})^m $, where $(\mathbb{R}/\mathbb{Z})$ is identified with the circle $S^1= \lbrace z \in \mathbb{C} \ : \ |z| = 1 \rbrace$.

{\theorem {\normalfont{(\cite{maximaltori})}} Any one-dimensional representation $\Phi : T \rightarrow \mathbb{C}$ is given by
$$
\mathbf{z}= (z_j)_j \in (S^1)^m \mapsto \mathbf{z}^{\lambda} = z^{\lambda_1}_1 z^{\lambda_2}_2 \ldots z^{\lambda_m}_m
$$
where $\lambda_1, \lambda_2, \ldots, \lambda_m$ are integers. 
}
\proof Let $\mathbf{t} = \mathbb{R}^m$ be the Lie algebra of $T$. The Lie algebra map $\phi : \mathbb{R}^m \rightarrow \mathbb{C}$ is as follows 
$$
(t_1,t_2,\ldots,t_m) \mapsto \sum_{j} \Tilde{\lambda}_j z_j
$$
for some scalars $\Tilde{\lambda}_j$. Therefore
$$
\Phi(\mathbf{z}) = \Phi(exp(t_1,t_2,\ldots,t_m)) = exp(\phi(t_1,t_2,\ldots,t_m)) = exp( \sum_{j} \Tilde{\lambda}_j z_j) \quad \quad \dag
$$
If we give $\mathbf{t}$ the standard basis of $\mathbb{R}^m$, i.e. $t_j = \delta_{kj}$, then $exp(\mathbf{t}) = 1$ by definition of the exponential map. This also forces $\Phi(exp(\mathbf{t})) = e^{\Tilde{\lambda}_k} = 1$ from here we conclude $\Tilde{\lambda}_j = 2 \pi i \lambda_j$ for an integer $\lambda_j$. Therefore, $\text{from}^{\dag}$
$$
\Phi(\mathbf{z}) = \prod_{j} e^{2\pi i \lambda_j t_j} = \prod_j (z)^{\lambda_j}_{j}
$$
by definition of $z_j$ above. 

\subsection{Differential Algebras}
Let $\omega$ be a $C^{\infty}$ $k-$form on the manifold $M$. The Lie derivative of a differential form is defined in a similar way to the Lie derivative of a vector field, but we use the pullback instead of the push forward to compare nearby values. 
\definition {\it{Lie derivative}} \cite{LOR08}. Let $X$ be a $C^{\infty}$ vector field on $M$, the Lie derivative by $X$ is a graded algebra of order zero. 
$$
\mathcal{L}_{X}: \Omega^{k}(M) \longrightarrow \Omega^{k}(M)
$$
$$
\mathcal{L}_{X} \left(\omega(Y)\right) \coloneqq \left( \mathcal{L}_{X} \omega \right)(Y) + \omega \left( \mathcal{L}_{X} Y \right) \in \Omega^{k}(M)
$$
{\theorem {\normalfont{(\cite{LOR20})}} Let $X$ be a $C^{\infty}$ vector field on $M$, then the following applies}

\begin{enumerate}[label=(\roman*)]
\item {\it{For $f \in C^{\infty}(M)$, $\mathcal{L}_{X}f = Xf$.  }}
\item {\it{For a $C^{\infty}$ vector field $Y$ on $M$, $\mathcal{L}_{X} Y = \left[X,Y\right]$. }}
\item {\it{The distributive property of the Lie derivative }
$$
\mathcal{L}_{X} \left( \omega \wedge \alpha \right) = \left( \mathcal{L}_{X} \omega \right) \wedge \alpha + \omega \wedge (\mathcal{L}_{X} \alpha)
$$
}
\item {\it{(The product formula) If $\omega \in \Omega^{k}(M)$, then for $Y_1,\ldots,Y_k \in \mathcal{X}(M)$}
$$
\mathcal{L}_{X} \left( \omega(Y_1, \ldots, Y_k) \right) = \left( \mathcal{L}_{X} \omega \right)(Y_1, \ldots, Y_k) + \sum_{i=1}^{k} \omega \left( Y_1, \ldots, \mathcal{L}_{X} Y_i, \ldots, Y_k \right)
$$
}
\end{enumerate}

{\corollary Suppose $G$ is a compact Lie group with Lie algebra $\mathbf{g}$, let $\rho_{a}$ be a representation of $G$ on $\Omega(M)$ $\forall a \in G$.
$$
\rho_{a} \circ \mathcal{L}_{X} \circ \rho_{a}^{-1} = \mathcal{L}_{\operatorname{Ad}_{a}X}
$$
where $\operatorname{Ad}$ is the adjoint representation of $G$ on $\mathbf{g}$.
}

\definition {\it{Interior multiplication}} \cite{NAK}. If $X$ is a $C^{\infty}$ vector field on a manifold $M$, the interior multiplication by $X$ is a graded algebra of degree $-1$
$$
\iota_{X}: \Omega^{k}(M) \longrightarrow \Omega^{k-1}(M)
$$
$$
\left( \iota_{X} \omega \right)_p \left( v_1, \ldots, v_{k-1} \right) \coloneqq \omega_p \left( X_p, v_1, \ldots, v_{k-1} \right)
$$
for $\omega \in \Omega^{k}(M)$, $p\in M$ and $v_1, \ldots, v_{k-1} \in T_pM$. The interior derivative is also called the {\it{contraction}}.

{\theorem {\normalfont{(\cite{LOR08})}} Let $X$ be a $C^{\infty}$ vector field on a manifold $M$ and $\Omega(M)$ the de Rham complex of $C^{\infty}$ forms on $M$. }

\begin{enumerate}[label=(\roman*)]
\item {\it{The contraction $\iota_X \omega$ is linear over $C^{\infty}$ functions in both arguments: for $f \in C^{\infty}$ and $\omega \in \Omega^{k}(M)$}
$$
\iota_{fX} \omega = f \iota_{X} \omega, \quad \quad \iota_{X}(f\omega)= f \iota_{X} \omega 
$$
}
\item {\it{The contraction $\iota_{X}: \Omega(M) \rightarrow \Omega(M)$ is an antiderivation of degree $-1$}
$$
\iota_{X} \left( \omega \wedge \alpha \right) = \left( \iota_{X} \omega \right) \wedge \alpha + (-1)^{\operatorname{deg}(\omega)} \omega \wedge \left( \iota_{X} \alpha \right)
$$
}
\item {\it{The composites of contractions}
$$\iota_{X} \circ \iota_{Y} = \left[X,Y\right] $$
}
\item {\it{Cartan's homotopy formula}
$$
\mathcal{L}_{X} = \iota_{X}d + d \iota_{X}
$$
}
\end{enumerate}

{\corollary Both the exterior derivative $d: \Omega^{k}(M) \rightarrow \Omega^{k+1}(M)$ and the interior multiplication $\iota_{X}: \Omega^{k}(M) \rightarrow \Omega^{k-1}(M)$ commute with the Lie derivative $\mathcal{L}_{X}: \Omega^{k}(M) \rightarrow \Omega^{k}(M)$ 
$$
d\mathcal{L}_{X} = \mathcal{L}_{X}d, \quad \quad \iota_{X}\mathcal{L}_{X} = \mathcal{L}_{X}\iota_{X} 
$$
} 

{\corollary Suppose $G$ is a compact Lie group with Lie algebra $\mathbf{g}$, let $\rho_{a}$ be a representation of $G$ on $\Omega(M)$ $\forall a \in G$.
$$
\rho_{a} \circ \iota_{X} \circ \rho_{a}^{-1} = \iota_{\operatorname{Ad}_{a}X}
$$
where $\operatorname{Ad}$ is the adjoint representation of $G$ on $\mathbf{g}$.
}
\remark Since the exterior derivative $d$ is independent of the chosen basis, this implies,
$$
\rho_{a} \circ d \circ \rho_{a}^{-1} = d
$$

\example Let $M = \mathbb{R}^n$, for all $i,j$
$$
\mathcal{L}_{\partial/\partial x_i} dx^j = d \mathcal{L}_{\partial/\partial x_i} x^j = d \left( \frac{\partial x^j}{\partial x_i} \right) = d \delta^{j}_{i} = 0 
$$
$$
\mathcal{L}_{\partial/\partial x_i} (f \ dx^j) =  \left( \frac{\partial f}{\partial x_i} \right)  x^j + f \ \mathcal{L}_{\partial/\partial x_i} dx^j = \left( \frac{\partial f}{\partial x_i} \right)  x^j
$$
$$
\iota_{\partial/\partial x_i} dx^j = dx^j \left( \frac{\partial }{\partial x_i} \right) = \left( \frac{\partial x^j}{\partial x_i} \right) =  \delta^{j}_{i}
$$
$$
\iota_{\partial/\partial x_i} (f\ dx^j) = f \ \iota_{\partial/\partial x_i} dx^j = f\delta^{j}_{i}
$$


\subsection{Hamiltonian Action and Moment Maps}

\definition \emph{Symplectic Manifold }\cite{symplctic}. A symplectic manifold $(M,\omega)$ is a smooth manifold equipped with a non degenerate two-form $\omega$ on $M$.
\remark In local coordinates, a general expression for the symplectic $2-$form on a $2n-$dimensional manifold $M$ can be written as
$$
\omega= \sum_{i=1}^{n} = dX_{i} \wedge dY_{i}
$$

\definition \emph{Symplectic volume }\cite{symplctic}. Since the top form on $M$ is a volume form, we define the symplectic volume as
$$
\frac{\omega^n}{n!} = dX_{1} \wedge dY_{1} \wedge \ldots \wedge dX_{n} \wedge dY_{n}
$$

Whenever we have an action of a Lie group $G$ on a symplectic manifold $(M, \omega) $ we also get the idea of a moment map.
\definition \emph{Moment map } \cite{holm2007act}. The moment map $\mu :  M\rightarrow \mathbf{g}^*$ which assigns to each pint $p\in M$ a map $\mathbf{t}\rightarrow \mathbb{R}, \ X \mapsto \phi^X(p)$ . A action for which a moment maps exists is called a Hamiltonian action; we may say a Lie group acts on a manifold in a Hamiltonian way. Furthermore, an action is said to be symplectic if it preserves the symplectic form: $g\omega = \omega \, \forall g\in G.$ 

\definition \emph{Fundamental vector field} \cite{holm2007act}. The fundamental vector field generated by an element of the Lie algebra $X \in \mathbf{g} $ is given assigning to each point $p\in M$ a vector
$$ \mathcal{X}_X(p) = \frac{d}{dt}[exp(tX)\cdot p] |_{t=0} \,.$$

\noindent This vector field is Hamiltonian if the contraction of the symplectic form by $\mathcal{X}_X$ is an exact form, 

$$ \iota_{X}\omega = d\mu^X \quad \quad (\dag)$$

\remark By virtue of the existence of \emph{maximal tori}, a Hamiltonian action is a symplectic torus action with moment map $\mu : M \rightarrow \mathbf{t}^{\vee}$, which also $\text{satisfies}^{\dag}$. \newline

 A useful technique when working with symplectic manifold with a $G$ action is symplectic reduction. 

{\theorem\label{symplectic_red}{\normalfont{(\cite{jose_lecture_3})}} Consider a Hamitonian $G$-action on the symplectic manifold $(M, \omega)$ with moment map $\mu.$ Let $\iota : \mu^{-1}(0) \hookrightarrow (M, \omega)$ be the inclusion map. Define the symplectic reduction of $M$ as the quotient $M_{red} =\mu^{-1}(0) / G $ and let $\pi : \mu^{-1}(0)  \rightarrow \mu^{-1}(0) / G$ be the natural projection map. Then $(M_{red}, \omega_{red})$ is a symplectic manifold with symplectic form $\omega_{red}$ characterised by $\iota^*\omega = \pi^*\omega_{red}$.} 

\subsection{Bundles, Connections and Curvature } 
The basic geometrical structure that we will be working with is that of a fibre bundle. They are foundational to the study of geometry and gauge theories is theoretical physics, as allow for a rigorous formulation of curvature, connections and covariant derivatives. The discovery that the curvatures of differential geometry are precisely the potentials in particle physics has led to a great deal of insight. 

\definition {\it Fibre Bundles:} \cite{Sch2017} A bundle is a triplet $(E, \pi, M)$ with $E$ called the total space, $M$ the base space and $\pi$ the projection map. $\pi : E \rightarrow M $ a continuous, surjective map. The fibre at at point $p\in M $ is the preimage $F_p = preim_{\pi}(\lbrace p \rbrace )$. 

\remark If $F$ is a manifold and $\forall p\in M \ preim_{\pi}(\lbrace p \rbrace) \cong F $ then $ (P, \pi, M)$ is a fibre bundle with typical fibre $F.$

\remark Of particular concern to this work will be vector bundles, where the total space is a vector space.

\definition {\it Principle $G$-Bundle:} \cite{Sch2017} Let G act freely to the right on $P$ and $\rho$ be the quotient map that sends $p\in P $ to its orbit $[p] \in P/ G$. Then $ (P, \pi, M)$ is a principle $G$-bundle if there is an isomorphism of bundles $(P, \pi, M) \cong (P, \rho, P/ G)$. In this case the fibres of $(P, \pi, M)$ are copies of $G$.

\remark A bundle is trivial if isomorphic to a product bundle $(M\times N, \pi ,M)$ where $\pi (m , n) = m$. The product bundle has fibres $N$, and one may think of taking the product of a manifold with its fibre as attaching a copy of the fibre at each point of $M$.

\definition \emph{section }\cite{NAK}. A section of a bundle $(P, \pi, M) $ is a map $\sigma : M \rightarrow P $ satisfying $$\pi \circ \sigma = id_M .$$

In the case of a vector bundle a section defined as vector field over the base space by assigning an element of the vector space to each point. The space of smooth sections for a vector bundle $(V, \pi, M)$ is denoted by $\Gamma (V)$.
  
\remark In gauge theories, such as Yang-Mills theory, one defines the fields and field strength on the base space through the construction of connections and curvature forms on the total space, which are then pulled back to the base space by a smooth section. 

{\definition {\it Vertical Subspace }\cite{NAK}. Let $(P, \pi, M)$ be a principal $G$-bundle and $p\in P.$ The vertical subspace at $p$ it the subspace of $T_pP$ tangent to the fibre $G.$ This is isomorphic to $\mathbf{g}$ at each point through the mapping of each element $X \in \mathbf{g}$ to its fundamental vector field.}

\definition {\it Connection One-Form }\cite{NAK}. A connection on $P$ is a unique direct sum decomposition $T_pP = V_pP \oplus H_pP$ which may be defined through a Lie algebra valued one-form $\omega \in \mathbf{g}\otimes T^*P $. $A$ acts as a projection of $T^*P$, meaning that $A(\mathcal{X}_X)= X$ and the pull-back of the right action of $G$ on $M$ is defined as $g^*\omega = g^{-1}\omega g$ for $g\in G$.

\definition \emph{The covariant derivative }\cite{diff_geom_symAgebra}. The covariant derivative $\nabla$ of a vector valued $k$-form $\alpha\in\Omega^k(P)\otimes V$ is defined for vector fields $X_1, ..., X_K \in \Gamma(TP)$ with horizontal projections $X_1^h, ...,X_k^h $ by $$ \nabla\alpha(X_1, .., X_k) := d\alpha(X_1^h, ...,X_k^h ) \,. $$

\definition \emph{Curvature }\cite{Sch2017}. The curvature $F_{\omega}$ of a connection $\omega$ is defined as  
$$F_{\omega} = d\omega + \omega \wedge \omega \,. $$ 

The second term is defined by its application to two vector fields $X, Y \in \Gamma(TP)$
$$(\omega \wedge \omega)(X, Y) = [\omega(X) , \omega(Y)] $$ 

\proposition The curvature two-forms $F_{\omega}$ are closed on $M$. 
$$ dF_{\omega} = d(dF_{\omega}) + d\omega \wedge \omega + \omega \wedge d\omega = 0$$

Usually in physics, one is interested in defining concepts locally on the base space, thus we use open coverings of the base space and we can proceed if out base space does not permit a global section.

\subsection{Introduction to Equivariant Cohomology}

Equivariant cohomology is a powerful topological invariant that captures essential information about group actions on spaces. Rather than being a property of the spaces themselves, it encodes intricate details about the space's topology, the isotropy groups associated with the action, and the arrangement of orbits, particularly focusing on the fixed points of the action. Pioneered by Borel and H. Cartan in the 1950s, equivariant cohomology has since found wide-ranging applications in diverse fields where symmetries of geometric objects are significant. The aim of these notes is twofold: to provide a gentle introduction to this elegant theory, approached through the lens of de Rham theory, and to survey both classical and contemporary applications that showcase its versatility and impact\cite{Goertsches_Zoller_2019}.
 
\definition {\it Universal Bundle} \cite{LOR08}. Let $G$ be a compact Lie group. The universal bundle $EG$ is a contractible space on which $G$ acts freely.

\definition {\it Classifying Space} \cite{Pestun_2017}. The space $BG = EG/G$ is called the classifying space of $EG$

{\theorem [\cite{Algeb_top}] For every topological group $G$, there exists a principal $G$-bundle bundle $EG \rightarrow BG$ where $EG$ is a contractible space. This leads to the concept that for every topological group $G$, there exists a unique classifying space $BG$ up to a homotopy quotient. }

\example Suppose $G=(\mathbb{Z},+)$, the universal bundle $EG = \mathbb{R}$ and the classifying space $BG=EG/ \mathbb{Z} = S^1$

\definition {\it{Borel Construction}} \cite{goertsches2019equivariant}. For any $G$-space $M$, the diagonal action of $G$ on $M \times EG$ is free. Giving rise to the \emph{homotopy quotient} $(M \times_{G} EG)$
$$
M_G = (M \times_{G} EG)  = (M \times EG)/G
$$

\definition {\it $G$-Equivariant Cohomology} \cite{LOR08}. For any $G$-space $M$, the $G$-equivariant cohomology of $M$ is the ordinary cohomology of the homotopy quotient $M_{G}$.
$$
H^{\bullet}_{G}(M;R) = H^{\bullet}(M_G;R) = H^{\bullet}((M \times EG)/G;R)
$$
{\corollary If $G$ acts freely on $M$ then the orbit space $M/G$ is a smooth manifold 
 {\normalfont{\cite{Pestun_2017}}}, thus}
$$
H^{\bullet}_{G}(M;R) = H^{\bullet}(M/G;R) 
$$

\remark The $G$-equivariant cohomology of a point is the ordinary cohomology of the classifying space.
$$
H^{\bullet}_{G}(pt;R)= H^{\bullet}((pt \times EG)/G;R) = H^{\bullet}(EG/G;R) = H^{\bullet}(BG;R)
$$

\example  For a compact torus $T= (S^1)^n$ we can take $ET=(S^{\infty})^n$, where $S^{\infty} \subseteq \mathbb{C}^{\infty}$ is the contractible infinite dimensional sphere. The Classifying space $BT=ET/T = (\mathbb{C}\mathbb{P}^{\infty})^n$ is $n$ copies of infinite complex projective space. The $T$-equivariant cohomology of a point is
\begin{align*}
H^{\bullet}_{T}(pt;\mathbb{Z})= H^{\bullet}(BG;\mathbb{Z})&= H^{\bullet}(\mathbb{C}\mathbb{P}^{\infty};\mathbb{Z}) \otimes \cdot\cdot\cdot \otimes H^{\bullet}(\mathbb{C}\mathbb{P}^{\infty};\mathbb{Z})\\
&= \mathbb{Z}[x_1]  \otimes \cdot\cdot\cdot \otimes \mathbb{Z}[x_n]\\
&\cong \mathbb{Z}[x_1,...,x_n]
\end{align*}
the polynomial ring in $n$ variables, where $deg(x_i) = 2$
{\theorem {\normalfont{(\cite{holm2007act})}} For any $G$-space $M$, we have the fibration
$M \hookrightarrow M_G \xrightarrow{p} BG$.
The projection $p$ induces $p^{*}: H^{\bullet}_{G}(pt;R) \rightarrow H^{\bullet}_{G}(M;R)$ making $H^{\bullet}_{G}(M;R)$ an $H^{\bullet}_{G}(pt;R)$ module. Natural maps in equivariant cohomology preserve this module structure}. 

{\theorem {\normalfont{(\cite{LOR08})}} Let $M, N$ be compact oriented $G$-manifolds, a diffeomorphism $f: N \rightarrow M$ induces the push-forward map $f_{*}: H^{\bullet}_{G}(N;\mathbb{Q}) \rightarrow H^{\bullet}_{G}(M;\mathbb{Q})$.}

\remark Consider the projection map $\pi: M \rightarrow pt$, the induced push-forward map $\pi_{*}: H^{\bullet}_{G}(M;\mathbb{Q}) \rightarrow H^{\bullet}_{G}(pt;\mathbb{Q})$ can be thought of as the equivariant integral.\newline

Let $\mathbf{g}= Lie(G)$ be the the complex Lie algebra of a compact Lie group $G$. Let $\mathbb{C}[\mathbf{g}]$ be the ring of complex valued polynomials on $\mathbf{g}$, and let $\mathbb{C}[\mathbf{g}]^G$ be the ring of $Ad_G$ invariant polynomials on $\mathbf{g}$.
$$
\mathbb{C}[\mathbf{g}]^G = \lbrace   f\in \mathbb{C}[\mathbf{g}] : f(Ad_{\mathbf{g}}x) = f(x) , \forall \mathbf{g} \in G, \forall x \in \mathbf{g} \rbrace \subseteq \mathbb{C}[\mathbf{g}]
$$
\definition {\it{Chern-Weil morphism}} \cite{Pestun_2017}. The characteristic classes $\mathbb{C}[\mathbf{g}]^G$ is isomorphic to the ordinary cohomology ring of the classifying space $BG$.
$$
H^{\bullet}_{G}(pt; \mathbb{C}) \cong  H^{\bullet}(BG; \mathbb{C}) \cong \mathbb{C}[\mathbf{g}]^G
$$
\example For the circle group, $G=S^1\cong U(1)$, $ES^1 = S^{\infty}$, $BS^1= \mathbb{C}\mathbb{P}^{\infty}$

$$
\mathbb{C}[\mathbf{g}]^G \cong H^{\bullet}(\mathbb{C}\mathbb{P}^{\infty}; \mathbb{C}) \cong \mathbb{C}[x]
$$
Where $x\in \mathbf{g^{\vee}}$ is a linear function on $\mathbf{g}= Lie(S^1)$ and $\mathbb{C}[x]$ is the polynomial ring with one generator $x$. Here $x \in H^{2}(\mathbb{C}\mathbb{P}^{\infty}; \mathbb{C})$ is, in fact, the first Chern class $c_1$ of the universal bundle up to a sign.\newline

{\theorem {\normalfont{(\cite{Chevalley})}} For a compact connected Lie group $G$ {\it{Chern-Weil morphism}} is reduced to the maximal torus $T\subset G$

$$
\mathbb{C}[\mathbf{g}]^G \cong \mathbb{C}[\mathbf{t}]^{W_G}
$$
Where $\mathbf{t}$ is the Lie algebra $\mathbf{t} = Lie(T)$ and $W_G$ is the Weyl group of G.}

\example Let $G=U(n)$ the Weyl group $W_{U(n)}$ is the permutation of group of $n$ eigenvalues $(x_1, \cdot \cdot \cdot, x_n)$. Therefore,

$$
H^{\bullet}(BU(n); \mathbb{C}) = \mathbb{C}[\mathbf{g}]^{U(n)} \cong \mathbb{C}[x_1, \cdot \cdot \cdot, x_n]^{W_{U(n)}} \cong \mathbb{C}[c_1, \cdot \cdot \cdot, c_n]
$$
Where $(c_1, \cdot \cdot \cdot, c_n)$ are symmetrical monomials called {\it{Chern classes}}
$$
c_k = (-1)^k \sum_{i_1 \leq \cdot\cdot\cdot \leq i_k} x_{i_1}\cdot\cdot\cdot x_{i_k}
$$


\subsection{Characteristic Classes and their Significance}
Characteristic classes are invariants used in algebraic topology to study and classify fibre bundles. They provide essential information about the topological properties of these bundles and help distinguish different bundles based on their algebraic structures. By associating characteristic classes with fibre bundles, we can gain insights into the fundamental properties and relationships of topological spaces. These classes play a crucial role in understanding the structure and behaviour of complex spaces and are invaluable tools in the classification of topological spaces.

\definition {\it{Characteristic classes }}\cite{Milnor}. A characteristic class is a map $c$ from a bundle $\xi = (E, \pi, M)$ into the ordinary cohomology of $M$, $c: \xi \rightarrow H^{\bullet}(M)$ which satisfies the naturaliy condition $c(\xi) = f^{*}c(\xi')$ for any bundle map $f: \xi \rightarrow \xi'$.

\definition {\it{Invariant polynomial }}\cite{Characteristic}. Let $M(n,\mathbb{C})$ be the space of $n \times n$ complex matrices. An invariant polynomial on $M(n, \mathbb{C})$ is a function
$$
P: M(n, \mathbb{C}) \rightarrow \mathbb{C}
$$
which is basis invariant $P(TAT^{-1})=P(A)$ for every nonsingular matrix $T$.
\remark The property of basis independence gives rise to the cyclicity $$P(AB)=P(BA)$$

\example Both the trace and determinant are invariant polynomials.\newline

Given a curvature $F_{\omega}$ ($\mathbf{gl}(n,\mathbb{C})$-valued 2-form on $M$) and any invariant polynomial $P$, we can define a differential form $P(F_{\omega})$ as follows. Consider an open cover of $M$ and in each open set, choose a local basis of sections $\lbrace\sigma_i\rbrace$. We may define the components $(F_{\omega})_{ij}$ of our curvature form in our chosen basis via
$$
F_{\omega}(\sigma_i) = \sum_{j} (F_{\omega})_{ij} \otimes \sigma_{j}
$$
Where each $(F_{\omega})_{ij}$ is a 2-form. The curvature form as a matrix $F_{\omega}=[(F_{\omega})_{ij}]$ whose entries lie in the commutative algebra of even-dimensional forms over $\mathbb{C}$, we can evaluate $P$ of $F_{\omega}$ (precisely due to the commutativity of all the elements of $F_{\omega}$). We may expand our polynomial as a power series,
$$
P = P_0 + P_1 + P_2 + \cdot\cdot\cdot
$$
where each $P_k$ is an invariant polynomial of degree $k$. The evaluation of $P(F_{\omega})$ will be will defined and $P_{k}(F_{\omega})=0$ for $2k > \operatorname{dim}M$ (since $P_{k}(F_{\omega})$ will be a differential form of degree $2k$)
\lemma{For any invariant polynomial $P$, the form $P(F_{\omega})$ is closed.}

{\corollary Since every invariant polynomial $P$ allows us to assign cohomology classes to curvatures on $\xi$:
$$
P: F_{\omega} \mapsto [P(F_{\omega})] \in H^{\bullet}(M;\mathbb{C})
$$
The cohomology class $ [P(F_{\omega})] \in H^{\bullet}(M;\mathbb{C})$ is independent of the connection $\omega$.
}
\definition{\it{Constructing curvature invariants }}\cite{charc_chern}. For any $(n \times n)$ matrix $A$, we define $\varepsilon_k(A)$ as the $k$-th elementary symmetric function of the eigenvalues of $A$:
$$
\varepsilon_{k}(A) = \sum_{1 \leq i_1 \leq \cdot\cdot\cdot \leq i_k \leq n} x_{i_1}\cdot\cdot\cdot x_{i_k}
$$
where $\lbrace x_i \rbrace$ are the eigenvalues of $A$. This implies that 
$$
\operatorname{det}(I+tA) = \sum_{k=0}^{\infty} t^k \varepsilon_k(A)
$$
Noting that $\varepsilon_{k}(A) = 0$ for $k > n$ so the sum is, in fact, finite. 

{\lemma Any invariant polynomial on $M(n,\mathbb{C})$ can be expressed as a polynomial function of $\varepsilon_{1} \cdot\cdot\cdot \varepsilon_{n}$.}

\example {\it{Chern classes}} $c_k$ are curvature invariant polynomials
$$
c_k(\xi) = \frac{\varepsilon_k(F_{\omega})}{(2\pi i)^k}
$$
\subsubsection*{The Euler Class}

For any oriented manifold, there is a privileged element of the top cohomology group (with coefficients in $\mathbb{Z}$) called the {\it{Euler class}}. For any bundle $\xi = (E,\pi, M)$ we will denote $E_0$ for the space obtained by removing the zero section from $E$. 

{\theorem {\normalfont{(Thom isomorphism \cite{Characteristic})}} Let $\xi =(E,\pi, M)$ be a real oriented $n$-plane bundle. The cohomology group $H^{i}(E,E_0;\mathbb{Z})$ is zero for $i < n$ and $H^{n}(E,E_0;\mathbb{Z})$ contains a unique cohomology class $u$ (the \textbf{Thom class}) such that the restrictions
$$
u|_{(F,F_0)} \in H^{n}(F,F_0;\mathbb{Z})
$$
are equal to the fundamental classes $u_F$ for each fibre $F$. The correspondence $y\mapsto y \cup u$ provides an isomorphism between $H^{k}(E;\mathbb{Z})$ and $H^{k+n}(E,E_0;\mathbb{Z})$ for each integer $k$. Therefore, $H^{\bullet}(E,E_0;\mathbb{Z})$ is a free $H^{\bullet}(E;\mathbb{Z})$-module on degree $n$ generator $u$.}

\definition {\it{Euler classes} }\cite{charc_chern}. The Euler class of an oriented (real) $n$-plane bundle $\xi =(E,\pi, M)$ is the cohomology class 
$$
e(\xi) \in H^{n}(M;\mathbb{Z})
$$
given by
$$
\pi^{*}e(\xi) = u|_E
$$
with the following properties

\proposition {\it{Naturality}}. For all $f: M \rightarrow M'$ covered by an orientation preserving bundle map $\xi \rightarrow \xi'$,
$$
e(\xi) = f^{*}e(\xi')
$$
\proposition {\it{Parity}}. If the orientation of $\xi$ is reversed, the Euler class changes sign. 

\example Let $\xi$ be a trivial $n$-plane bundle $M \times \mathbb{R}^n \xrightarrow{\pi} M$. Consider the bundle map that takes $\xi$ to the trivial $n$-plane bundle over a point. 
$$
f: M \times \mathbb{R}^n \rightarrow \lbrace 0\rbrace \times \mathbb{R}^n
$$
Here $M$ is mapped into $\lbrace 0\rbrace$, this mapping is orientation preserving by construction. 
$$
e(\xi)=(\pi \circ f)^{*} e(\xi')= (\pi \circ f)^{*}0 = 0
$$
\proposition The Euler class of a Whitney sum or a Cartesian product of bundles
$$
e(\xi \oplus \xi') = e(\xi)  e(\xi') \quad\quad\quad  e(\xi \times \xi') = e(\xi) \times e(\xi')
$$

\lemma If the fibre dimension $n$ of $\xi$ is odd, $e(\xi)+e(\xi)=0$

\lemma Let $\xi = (E, \pi, M)$ be a real oriented vector bundle with a nowhere zero section, $e(\xi) = 0$

\subsubsection*{The Chern Class}

\definition {\it{Chern classes }}\cite{Milnor}. The Chern classes of a complex $n$-dimensional bundle $\zeta$ are the generators of the cohomology ring $c_k(\zeta) \in H^{2k}(M;\mathbb{Z})$\newline

\noindent For $\mathbf{k>n}$ we set $c_k(\zeta) = 0$. The top Chern class $c_n(\zeta)$ is equal to the Euler class of the underlying real bundle. 
$$
c_n(\zeta) = e(\zeta_{\mathbb{R}}) \in H^{2k}(M;\mathbb{Z})
$$

\noindent For $\mathbf{k<n}$, we consider the map $(\pi^{*}_0)^{-1}: H^{2k}(E_0) \rightarrow H^{2k}(M)$, we set 
$$
c_k(\zeta) = (\pi^{*}_0)^{-1} c_k(\zeta)
$$
For example, this sets 
$$
c_{n-1}(\zeta) = (\pi^{*}_0)^{-1} e(\zeta_{0,\mathbb{R}})
$$
In particular, this gives $c_0(\zeta) = 1$
{\corollary {\normalfont{Total Chern class}} of a complex $n$-dimensional bundle $\zeta$ over $M$ is defined as}{ \cite{Characteristic}.}
$$
c(\zeta) = \sum_{k=1}^{n} c_k(\zeta) \in \bigoplus_{k}^{n} H^{k}(M;\mathbb{Z})
$$
with the following properties,
\proposition {\it{Naturality}}. If $f: M \rightarrow M'$ is covered by a bundle map from a complex $n$-plane bundle $\zeta$ to a complex $n$-plane bundle $\zeta'$, then 
$$
c(\zeta) = f^{*} c(\zeta')
$$
\proposition {\it{Sum formula}}. Let $\zeta$ and $\eta$ be two complex vector bundles over a common paracompact base space $M$, 
$$
c(\zeta \oplus \eta) = c(\zeta)  c(\eta)
$$

\remark If $\epsilon^k$ is the trivial complex $k$-plane bundle over $M$, then 
$$c(\zeta \oplus \epsilon^k)= c(\zeta)$$.

{\theorem{\normalfont{(Splitting Principle \cite{Milnor})}}. Let $\zeta = (E, \pi, M)$ be a complex $n$-plane bundle over $M$, there exits a morphism $f: M' \rightarrow M$ such that the induced pull-back $f^{*}: H^{\bullet}(M) \rightarrow H^{\bullet}(M')$ is injective, and
$$
f^{*}(\zeta) \cong L_1 \oplus \cdot\cdot\cdot \oplus L_n
$$
where $L_i$ are complex line bundles over $M'$}

{\corollary The first Chern class is the unique (modulo scalar multiplication) nontrivial characteristic class of a complex line bundle}{ \cite{charc_chern}.}

\remark If $\zeta$ is a complex vector bundle of rank $n$, then $$c_{n}(\zeta) = e(\zeta_{\mathbb{R}})$$

\proposition If $\zeta = L_1 \oplus \cdot\cdot\cdot \oplus L_n$ (direct sum of line bundles) then by the splitting principle, 
$$e(\zeta) = c_1(L_1)\cdot\cdot\cdot c_1(L_n)$$

\example Consider a closed, oriented Riemannian two-manifold $M$. In two dimensions, the Levi-Civita curvature two-form on $M$ is 
$$
(F_{\omega})_{ij} = K \theta_i \wedge \theta_j
$$
where $K$ is the Gaussian curvature of $M$.\newline
In a neighbourhood $U$ of $M$ we can introduce geodesic coordinates $(x,y)$ in which the \textbf{first fundamental form} (the Euclidean metric) takes the form
$$
\mathbf{I} = dx \otimes dx + g(x,y)^2 dy \otimes dy
$$
Setting 
$$
\theta_1 = dx \quad\quad\quad \theta_2 = g(x,y) dy
$$
gives local orthonormal basis over $U$. We can consider the first structure equation
$$
0 = d\theta_1 = \omega_{12} \wedge \theta_2 = g(x,y) \omega_{12} \wedge dy
$$
$$
\frac{\partial g(x,y)}{\partial x} dx \wedge dy = d\theta_2 = \omega_{21} \wedge \theta_1 = - \omega_{12} \wedge \theta_1 = - \omega_{12} \wedge dx
$$
Therefore, the connection one-form is given by
$$
\omega_{12} = \frac{\partial g(x,y)}{\partial x} dy
$$
which implies
$$
(F_{\omega})_{12} = \frac{\partial^2 g}{\partial x^2} dx \wedge dy = \frac{1}{g} \frac{\partial^2 g}{\partial x^2} d\theta_1 \wedge d\theta_2 = - \frac{1}{g} \frac{\partial^2 g}{\partial x^2} \theta_1 \wedge \theta_2 
$$
Thus, we must have 
$$
K = - \frac{1}{g} \frac{\partial^2 g}{\partial x^2}
$$

{\theorem {\normalfont{(Gauss-Bonnet \cite{MartensGeo})}}. Let $M$ be a closed, oriented Riemannian two-manifold, then}
$$
\int_{M} K dA = 2\pi \chi(M)
$$
{\it{Proof}}: The orientation of $M$ gives rise to an orientation to each of the fibres of $TM$ which we can view as an oriented two-plane bundle with a Euclidean metric. Let's consider an arbitrary oriented, real two-plane bundle $\xi$, it has a canonical complex structure $J$ which rotates a vector "counterclockwise" by $\pi/2$ i.e. in terms of any oriented local orthonormal basis $\lbrace \sigma_i \rbrace$, $J \sigma_1 = \sigma_2$. We may canonically identify $\xi$ with a complex line bundle $\zeta$ which inherits a connection obeying  
$$
\nabla(\sigma_1) = \omega_{12} \otimes \sigma_2 = \omega_{12} \otimes (i\sigma_1) = i \omega_{12} \otimes \sigma_1
$$
so that 
$$
\nabla(i\sigma_1) = - \omega_{12} \otimes \sigma_1
$$
therefore, the connection one-form on $\zeta$ is given by
$$
\omega_{\zeta} = i \omega_{12}
$$
thus the curvature two-form on $\zeta$ is given by
$$
(F_{\omega})_{\zeta} = i (F_{\omega})_{12}
$$
Using the invariant polynomial $\varepsilon_1 = \operatorname{tr}$, we obtain the cohomology class
$$
[\operatorname{tr}(i (F_{\omega})_{12})] = [i (F_{\omega})_{12}] \in H^2 (M;\mathbb{C})
$$
This represents a characteristic class for $\zeta$, here we can use the existence and uniqueness of the Chern class
$$
i (F_{\omega})_{12} = \alpha c_1(\zeta) = \alpha e(\xi)
$$
for some $\alpha \in \mathbb{C}$. To evaluate $\alpha$, it suffices to check the case of $TM$, since $i(F_{\omega})_{12} = \alpha e(TM)$. Integrating both sides we must have 
$$
i \int_{M} K dA = \int_{M} i(F_{\omega})_{12} = \alpha \int_{M}  e(TM) 
$$
Since $\int_{M}  e(TM) = \chi(M)$ is the fundamental reason we name the top Chern class {\it{Euler class}}

$$
i \int_{M} K dA = \alpha \chi(M)
$$
We can evaluate both sides for a sphere,  $K=1$ and $\chi(S^2) = 2$, to get
$$
4\pi i = 2 \alpha \quad \implies \quad \alpha = 2\pi i
$$
Hence, 
$$
\int_{M} K dA = 2\pi \chi(M) \quad \square
$$

% ...

\newpage

\section{Chapter 3: Equivariant Cohomology and Localization}
\subsection{Cartan Model of Equivariant Cohomology}

We have introduced the de Rham cohomology of a manifold, which is a differential graded algebra. If the manifold has a Lie group action defined on it then the de Rham complex will have additional algebraic structure. The algebraic properties of the equivariant differential forms are expressed in a $\mathbf{g}$-differential graded algebra. This is constructed through the Weil model and the Cartan model, which are isomorphic; however, the later gives a simpler formulation so is more widely used.

\definition \emph{Cartan algebra }\cite{Meinrenken_Equivariant_Cohomology}. Let $\Omega(M)$ be the de Rham graded algebra of a $G$-maniold $M$.
$$
\Omega_{G}(M) = \left( S(\mathbf{g}^{\vee}) \otimes \Omega(M)  \right)^G 
$$
\noindent $S(\mathbf{g}^{\vee})$ is the symmetric algebra of $\mathbf{g}^{\vee}$, where $\mathbf{g}^{\vee}$ is the dual of $\mathbf{g} \in G$
$$
S(\mathbf{g}^{\vee}) \coloneqq \lbrace f: \mathbf{g} \rightarrow \mathbb{R} \ | \   f \text{ is a polynomial} \rbrace 
$$
\definition \emph{Cartan graded algebra }\cite{LOR20}. Let $\Omega_{G}(M)$ be the Cartan algebra,
$$
\Omega_{G}(M) \coloneqq \bigoplus_{n=1}^{\operatorname{dim}(M)} \Omega_{G}^{n}(M)
$$
$$
\text{where } \quad \Omega_{G}^{n}(M) = \bigoplus_{n=2i+j} \left( S^{i}(\mathbf{g}^{\vee}) \otimes \Omega^{j}(M)  \right)^G 
$$

\remark If the Lie group is abelian, e.g. a torus $T$, then $Ad : T \lcirclearrowright\mathbf{t}$ is trivial thus $$ \left( S(\mathbf{g}^{\vee})\otimes\Omega(M) \right)^T = S(\mathbf{g}^{\vee})\otimes\left( \Omega(M) \right)^T \,.$$

\definition \emph{Cartan differential operator} \cite{Meinrenken_Equivariant_Cohomology}. Let $\lbrace \xi_i \rbrace$ be a basis for $\mathbf{g}$ and let $\lbrace \xi^i \rbrace$ be the dual basis. The Cartan differential operator is defined as 
$$
d_{G}: \Omega_{G}^{n}(M) \longrightarrow \Omega_{G}^{n+1}(M)
$$
$$
\text{by } \quad d_{G} := 1\otimes d - \sum_i \left( \xi^i \otimes\iota_{\xi_i} \right) \,.
$$

{\corollary The action of Cartan differential operator $d_{G}$ on an equivariant differential n-form $\omega \in \Omega^{n}_{G}(M)$ has degree $1$ in the Cartan graded algebra and gives rise to the Cartan Complex
$$
\xymatrix{
    {\Omega_{G}^{0}(M) \coloneqq \left( S^{0}(\mathbf{g}^{\vee}) \otimes \Omega^{0}(M)  \right)^G}  \ar@{->}[d]_{d_{G}} \\
    {\Omega_{G}^{1}(M) \coloneqq \left( S^{0}(\mathbf{g}^{\vee}) \otimes \Omega^{1}(M)  \right)^G} \ar@{->}[d]_{d_{G}} \\
    {\Omega_{G}^{2}(M) \coloneqq \left( S^{1}(\mathbf{g}^{\vee}) \otimes \Omega^{0}(M)  \right)^G \oplus \left( S^{0}(\mathbf{g}^{\vee}) \otimes \Omega^{2}(M)  \right)^G} \ar@{->}[d]_{d_{G}}|{\vdots\vdots\vdots} \\
    {\Omega_{G}^{n}(M) \coloneqq {\bigoplus\limits_{\substack{n=2i+j}} \left( S^{i}(\mathbf{g}^{\vee}) \otimes \Omega^{j}(M)  \right)^G}}  \\ 
}
$$
}

\proposition Analogous to de Rham operator, $d_G^2 = 0 $ on $\left( S(\mathbf{g}^{\vee})\otimes\Omega(M) \right)^G$.
\proof This proof uses the following identities:
$$
 d^2 = 0 \,,\quad \lbrace d, \iota_{\xi_i}\rbrace = \mathcal{L}_{\xi_i} \,,\quad \iota_{\xi_i} \circ \iota_{\xi_j} = 0 \, \quad \forall \xi_i \in \mathbf{g} \,.
$$
\begin{align*}
    d_G^2 &= 1\otimes d^2 -\sum_i \left( \xi^i  \otimes \lbrace d, \iota_{\xi_i} \rbrace \right) + \sum_i \sum_j \left( \xi^i \xi^j \otimes \iota_{\xi_i} \circ \iota_{\xi_j} \right) \\
    &=  - \sum_i \left( \xi^i \otimes \mathcal{L}_{\xi_i} \right) \\
    &= 0 \quad \quad \square
\end{align*}

{\lemma Cartan differential operator obeys the super Leibniz rule
$$
d_{G} \left(\alpha \wedge \omega\right) = \left(d_{G}\alpha\right) \wedge \omega + (-1)^{\operatorname{deg}\alpha} \alpha \wedge d_{G}\omega
$$
}

\proposition Cartan differential operator $d_{G}$ is basis independent,
\proof Let $\mathbf{g}$ be spanned by $\lbrace \xi_{i} \rbrace$ and $\lbrace \zeta_{j} \rbrace$, there exists transformation functions such that,
$$
\xi_i = U_{i}^{j} \zeta_{j}, \quad \quad \xi^i =  \zeta^{j}U^{i}_{j} ,
$$
substituting to the Cartan differential operator (using the summation convention)
\begin{align*}
d_{G} &= 1\otimes d -   \xi^i \otimes\iota_{\xi_i} \\
&=  1\otimes d -    \zeta^{j} U^{i}_{j} \otimes\iota_{U_{i}^{j} \zeta_{j}} \\
&= 1\otimes d -  \zeta^{j}  U^{i}_{j}  \otimes U_{i}^{j} \iota_{ \zeta_{j}} \\
&= 1\otimes d -  \zeta^{j}  \otimes  \iota_{ \zeta_{j}} \quad \quad \square
\end{align*}

\definition \emph{Cartan $G-$equivariant cohomology} \cite{Meinrenken_Equivariant_Cohomology}. In line with the initial definition of the cohomology ring, the $n$th $G-$equivariant cohomology ring is a measure of failure for an equivariant $n-$form to be exact.

$$
H^{n}_{G}(M) \coloneqq \frac{\operatorname{ker}d_{G}: \Omega^{n}_{G}\rightarrow \Omega^{n+1}_{G}}{\operatorname{im}d_{G}: \Omega^{n-1}_{G}\rightarrow \Omega^{n}_{G}}
$$
Equivalently, Cartan $G-$equivariant cohomology could be defined as,

$$ 
H^{\bullet}_{G}(M) \coloneqq Z^{\bullet}_{G}(M) / B^{\bullet}_{G}(M) =  \lbrace [\omega] : \omega\sim\mu \iff \omega - \mu \in B^{\bullet}_{G}(M) \rbrace \,.
$$
Where $\lbrace Z^{\bullet}_{G}(M), B^{\bullet}_{G}(M)\rbrace$ denote the sets of equivariant $\lbrace$closed, exact$\rbrace$ forms on $M$ respectively.

{\theorem {\normalfont{Cartan theorem \cite{Pes17}}}. $H^{\bullet}_{G}(M)$ is isomorphic to the cohomology of differential complex $H^{\bullet}\left( \Omega^{\bullet}_{G}(M), \ d_G \right) $.}
$$
H^{\bullet}\left( (M \times EG)/G, \ d \right) \ \scalebox{1.2}{$\cong$} \ H^{\bullet}\left( \left( S(\mathbf{g}^{\vee}) \otimes \Omega(M)  \right)^G,  d_{G}  \right)
$$

{\lemma {\normalfont{(\cite{Chevalley})}} In the case where $G$ could be reduced to a maximal torus, Chern-Weyl morphism could become:
$$
S(\mathbf{g}^{\vee})^G \cong S(\mathbf{t}^{\vee})^{W_G}
$$
}
{\lemma {\normalfont{(\cite{diff_geom_symAgebra})}} The fact that $\mathbf{g} \cong \mathbf{g}^{\vee}$ implies the following isomorphism
$$
S(\mathbf{g}^{\vee}) \cong \mathbb{C}\left[ \xi_1, \xi_2, \ldots, \xi_n \right]
$$
}

\example The equivariant cohomology of a point in the Cartan model
$$
H_G^{\bullet}(pt) = S(\mathbf{g}^{\vee})^G  \cong \mathbb{C}[\mathbf{g}]^G
$$
Therefore, the equivariant cohomology of a point in the Cartan model is compatible with that in de Rham model of equivariant cohomology shown earlier. \newline

Consider a symplectic manifold $(M,\omega)$ equipped with a Hamiltonian action of a compact group $G$ with a moment map $\mu$, we define 
$$
\overline{\omega}(X) = \omega + \mu \ \in \Omega^{2}_{G}(M) 
$$

{\proposition The equivariant Symplectic two-form is closed in the Cartan complex 
$$
d_{G} \overline{\omega}(X) = 0
$$}
\proof Fist, we decompose $\overline{\omega} = \omega + \mu$ to its principle components in the Cartan graded algebra. 
$$
\omega \in \left( S^{0}(\mathbf{g}^{\vee}) \otimes \Omega^{2}(M)  \right)^G, \quad \quad \mu \in \left( S^{1}(\mathbf{g}^{\vee}) \otimes \Omega^{0}(M)  \right)^G 
$$
Now let $\lbrace \xi_i \rbrace$ be a basis for $\mathbf{g}$ and let $\lbrace \xi^i \rbrace$ be the dual basis.

$$
\omega \rightarrow (1 \otimes \omega), \quad \quad \mu \rightarrow \sum_{i}(\xi^{i} \otimes \mu)
$$

\begin{align*}
d_{G} \overline{\omega}  &= (1 \otimes d)(1 \otimes \omega + \sum_{i}\xi^{i} \otimes \mu ) - (\sum_{i} \xi^{i} \otimes \iota_{\xi_{i}}) (1 \otimes \omega + \sum_{i}\xi^{i} \otimes \mu ) \\
&= \sum_{i}\xi^{i} \otimes d\mu \ - \ \sum_{i} \xi^{i} \otimes \iota_{\xi_{i}} \omega
\end{align*}
Since $\iota_{\xi_{i}} \omega = d\mu$ by definition of the Hamiltonian group action.
$$
\implies \  d_{G} \overline{\omega} = 0 \quad \quad \square
$$
Therefore $\left[\overline{\omega}\right] \in H^{2}_{G}(M)$

{\theorem {\normalfont{(\cite{goertsches2019equivariant})}} Let $M$ be equipped with a Hamiltonian $T-$action, and let $F$ be a component of the fixed points set on $M$. We have a natural inclusion map $\iota_F: F \rightarrow M$. For $\overline{\alpha} \in H^{\bullet}_{T}(M)$, the push-forward of such a form is 
$$
\iota^{*}_F \overline{\alpha} \in H^{\bullet}_{T}(F) \cong H^{\bullet}(F) \otimes H^{\bullet}_{T}(pt)
$$
}
{\lemma Let $(M,\omega)$ be a symplectic manifold with a Hamiltonian $T-$action. If $F \subset M$ is a component of the fixed point set then \label{symplectic_push} }
\begin{align*}
\iota^{*}_{F}\overline{\omega} &\in \Omega^{2}_{T}(F) \\
&= \omega|_{F} + \mu(F)
\end{align*}
\subsection{Equivariant Characteristic Classes}

\definition \emph{Equivariant Chern classes} \cite{Equi_Characteristic}. Let $G$ be a compact Lie group, and $\zeta=(E, \pi, M)$ be a complex vector bundle on manifold $M$.equipped with a Hamiltonian action. The Equivariant Chern classes $c^{G}_{k}(\zeta)$ are given by
$$
c^{G}_{k}(\zeta) \coloneqq c_{k}\left((\zeta \times_{G} EG) \longrightarrow (M \times_{G} EG)  \right)
$$
\definition \emph{Equivariant Euler classes} \cite{berlin_equivariant_char}. Since the Euler classes are the top Chern classes, therefore, the equivariant Euler classes are defined as
$$
e^{G}(\zeta) \coloneqq e \left((\zeta \times_{G} EG) \longrightarrow (M \times_{G} EG)  \right)
$$

\definition \emph{Equivariant characteristic classes in the Cartan Model} \cite{jeffery_QFT}. Suppose $\zeta$ is a complex vector bundle on $M$ equipped with the action of a compact Lie group $G$. Let $\omega$ be a connection on $\zeta$ compatible with the action of $G$. Define the moment of $\zeta$, $\mu \in \Gamma(\operatorname{End}(\zeta) \otimes \mathbf{g}^{\vee})$ as follows,
$$
\mathcal{L}_{X} \sigma - \omega_{X} \sigma = \mu(X) \sigma \quad \quad (\dag)
$$
for $\sigma \in \Gamma(\zeta)$ sections on the complex vector bundle $\zeta$, and $X \in \mathbf{g}$ basis of the Lie algebra of $G$. The action of $G$ allows us to define the Lie derivative $\mathcal{L}_{X} \sigma$ of the section $\sigma \in \Gamma(\zeta)$, and the formula $\text{above}^{\dag}$ defines $\mu$ as a zero order operator.   

\proposition The equivariant Chern class $c^{G}_{k}(\zeta)$ in the Cartan model are given by
$$
c^{G}_{k}(\zeta, X) = \left[ \varepsilon_{k}(F_{\omega} + \mu(X)) \right]
$$
where $F_{\omega} \in \Gamma(\operatorname{End}(\zeta) \otimes \Omega^{2}(M))$ is the curvature of $\omega$, and $\varepsilon_{k}$ is the $k-$th elementary symmetric polynomials.\newline

Now we examine the case where we have a fixed point set $F \subset M$

\definition\label{circ-action-weight} \emph{Weight of the action} \cite{jeffery_QFT}. A principal circle bundle over a point $pt$ equipped with the action of a torus $T$ is $P=S^1$ equipped with a weight $\beta \in \operatorname{Hom}(T,U(1))$. Thus, $T$ acts on $S^1$ by 
$$
t \in T: \ z \in S^1 \mapsto \beta(t)z
$$
\example  A principal circle bundle $S^1$ over a point $pt$ equipped with the action of a circle $U(1)$ has the weight $\beta(t) = t^m$.
$$
t \in U(1): \ z \in S^1 \mapsto t^{m} z \quad \text{for } m \in \mathbb{Z}
$$

\remark\label{first_cher} Since $\iota^{*}_{F}\overline{\omega} = \omega|_{F} + \mu(F)$ where $F \subset M$ is the fixed point set, for any $f\in F$, $\zeta|_f$ is a copy of $S^1$ on which $T$ acts using a weight $\beta_{F} \in \operatorname{Hom}(T, U(1))$. Here $\operatorname{Lie}(\beta_{F}) \in \operatorname{Hom}(\mathbf{t}, U(1))$. The equivariant first Chern class is
$$
c^{T}_{1}(\zeta|_{F}, X) = c_1(\zeta|_{F}, X)| - \operatorname{Lie}(\beta_{F}) = \left[\omega|_{F}\right]| - \operatorname{Lie}(\beta_{F})
$$
Therefore, at fixed points of the action, the value of the moment map is a weigh 
$$\mu(F) = -\operatorname{Lie}(\beta_{F})$$

\example Let $\zeta=(L,\pi,M)$ be a complex line bundle equipped with a Hamiltonian $T-$action, and $F \subset M$ is a fixed point set on $M$. Suppose $T$ acts on the fibres of $\zeta|_{F}$ with weights $\beta_{F} \in \mathbf{t}^{\vee}$ i.e. $exp(X) \in T: \ z \in \zeta|_{F} \mapsto e^{i \beta_{F}(X)}z$. Thus,

$$
e^{T}(\zeta|_{F},X) = c_{1}(\zeta|_{F}) + \beta_{F}(X)
$$
\remark Let $\zeta=(E,\pi,M)$ be a complex vector bundle equipped with a Hamiltonian $T-$action, and $F \subset M$ is a fixed point set on $M$. Let $\nu_F$ be the normal bundle to $F \subset M$, the $T$ acts on $\nu_F$. Using the splitting preiciple, $\nu_F$ decomposes $T-$equivariantly as a sum of line bundles $\nu_{F,j}$ on each of which $T$ acts with weight $\beta_{F,j} \in \mathbf{t}^{\vee}$. Therefore, the Euler class of the normal bundle $\nu_{F}$ is
$$
e^{T}(\nu_F,X) = \prod_{j} \left( c_1(\nu_{F,j}) + \beta_{F,j}(X)\right)
$$
\remark The weight $\beta_{F,j} \neq 0$ for any $j$ since $\nu_{F,j}$ is normal to the fixed point set. Thus
$$
e^{T}_{0}(\nu_F,X) = \prod_{j} \beta_{F,j}(X)
$$
which implies 
$$
e^{T}(\nu_F,X) = e^{T}_{0}(\nu_F,X)  \prod_{j} \left( 1+ \frac{c_1(\nu_{F,j})}{\beta_{F,j}(X)}\right)
$$
since ${c_1(\nu_{F,j})}/{\beta_{F,j}(X)}$ is nilpotent \cite{AtiBot84}, we can define the inverse of the Euler class
$$
\frac{1}{e^{T}(\nu_F,X)} = \frac{1}{e^{T}_{0}(\nu_F,X)} \sum_{k=0}^{\infty} (-1)^{k} \left(\frac{c_1(\nu_{F,j})}{\beta_{F,j}(X)}\right)^k
$$
only finite number of terms contribute to this sum.

\example The action of $T= U(1)$ on manifold $M$ with isolated fixed points $F\subset M$.

Suppose the normal bundle $\nu_{F} = T_FM$ at each point $F$, using the splitting principle the normal bundle decomposes as a direct sum $\nu_{F} \cong \bigoplus_{j=1}^{N} \nu_{F,j}$ where each $\nu_{F,j} \cong \mathbb{C}$ and $M$ acts with multiplicity $m_{F,j}$ on $\nu_{F,j}$ (a non zero integer)
$$
t \in U(1): \ z_j \in \nu_{F,j} \mapsto t^{m_{F,j}} z_j 
$$
therefore, the equivariant Euler class is 
$$
e^{T}(\nu_F,X) = \left(\prod_{j} m_{F,j} \right) X^{N}
$$

\subsection{Localization Techniques in Equivariant Cohomology} 
We will now introduce the central theorem of this report: the localization theorem. This had an early use in the proving the exactness of the statonary phase approximation and has now found a wide range of application in various forms, in particular to instanton counting. 

{\theorem\label{localisation_form} {\normalfont{(Atiyah-Bott Localization \cite{AtiBot84})}}. Let $(M, \omega)$ be a compact symplectic manifold with Hemiltonian $T$-action and moment map $\mu.$ Let $F \subset M$ be the fixed point set with elements $f \in F$, and $\nu_F$ be the normal bundle at $F$. }
$$\int_M \overline{\omega}(X)  = \sum_{f\in F } \int_F \frac{\iota^{*}_{F}\overline{\omega}(X)}{e^{T}(\nu_F,X)} ,. $$
When the fixed points are isolated, the integral reduces to evaluation at point.

{\lemma\label{power_lem} Let $M$ be a $2n-$dimensional manifold equipped with a Hamiltonian $T-$action. Assuming isolated fixed points, the localisation formula implies that for all integers $k \geq 0$,
$$
\int_M  \left(\omega + \mu(X) \right)^{k}  = \sum_{f\in F } \frac{\mu|_{f}(X)^k}{e^{T}(\nu_F,X)}
$$
where $\mu|_{f}$ denotes the moment map evaluated at point $f \in F$. For $k=n$, both sides are independent of $X$. Of course, the integral vanishes for $k<n$.
}
{\theorem \emph{(Duistermaat-Heckman Localization \cite{DuiHec82})} Let $(M, \omega)$ be a symplectic manifold of $\operatorname{dim}_{\mathbb{R}}(M) = 2n$ equipped with a Hamiltonian $T-$action whose moment map is $\mu.$ Assume that the fix points $F \subset M$ are discrete,
$$
\int_M \frac{(\omega)^n}{n!} \operatorname{exp}(\mu(X)) = \sum_{f\in F } \frac{\operatorname{exp}(\mu(X)) }{e^{T}(\nu_F,X)} \,.
$$}

\example A very simple demonstration of the localization theorem is the evaluation of the Gaussian integral. Let $T = U(1) = S^1$ act on the complex plane $\mathbb{C}$ and $t\in\mathbb{R}.$ The fixed point of the action is the origin and the weight is $t$, thus
$$
\int_{\mathbb{C}}e^{-t || z||^2}dz\wedge d\overline{z} = \frac{1}{t} \,.
$$
This is true whenever the integral converges, which occurs when $t>0.$\newline
{\theorem{{\normalfont{(Atiyah-Bott-Kirwan map \cite{jeffery_QFT})}} Suppose $(M,\omega)$ is a symplectic manifold equipped with a Hamiltonian action of a compact Lie group $G$ with a moment map $\mu$. The natural inclusion map $\iota: \mu^{-1}(0) \hookrightarrow M$ gives rise to a cohomology map
$$
\kappa : H^{\bullet}_{G}(M) \rightarrow H^{\bullet}_{G}(\mu^{-1}(0))
$$
Since the orbit space $\mu^{-1}(0)/G$ is an ordinary manifold we called $M_{red}$,
$$
\kappa : H^{\bullet}_{G}(M) \rightarrow H^{\bullet}(M_{red})
$$
}}
\remark Although the inclusion map $\iota$ is injective, the cohomology map $\kappa$ is surjective. This theorem has powerful applications in the ADHM construction and moduli spaces of instantons. \newline

\example Suppose $M = \mathbb{CP}^n$ be the $n-$dimensional complex projective space ($\operatorname{dim}_{\mathbb{R}M} = 2n$). Let $c^{T}_1(\mathcal{O}(1))$ be the $T-$equivariant first Chern class of the hyperplane bundle over $\mathbb{CP}^n$. Calculate the following integral using the localisation theorem.
$$
\int_{\mathbb{CP}^n} c^{T}_1(\mathcal{O}(1))^n
$$
The fixed points of $\mathbb{CP}^n$ under the action of the $n$-torus are the discrete set of points with homogeneous coordinates 
$$F \coloneqq(\mathbb{CP}^n)^T = \underbrace{\lbrace (1: 0: ... : 0), (0: 1: ... : 0), ... \,, (0: ... : 0: 1) \rbrace \,}_{(n+1) \text{ times}}.$$
The $T-$equivariant first Chern class $c^{T}_1(\mathcal{O}(1))$ is a generator of $H_{T}^2(\mathcal{O}(1); \mathbb{Z})$. Let $\omega$ be a symplectic form on $\mathbb{CP}^n$ with moment map $\mu(X),$ for $t\in \mathbf{t}$. By \emph{remark} \ref{first_cher},
$$
c^{T}_1(\mathcal{O}(1)) = [\omega] + \mu(t)
$$
Using the localisation formula \emph{theorem} \ref{localisation_form} and \emph{lemma} \ref{power_lem}
$$
    \int_{\mathbb{CP}^n}([\omega] +\mu(t))^n = \sum_{f\subset F } \frac{\mu|_{f}(t)^n}{e^{T}(\nu_F,t)}
$$
The left hand side is an element of the equivariant cohomology of a point $H_{T}^{\bullet}(pt) = \mathbb{C}[t_0, ... , t_n] $ which is isomorphic to the ring  of polynomials in $n+1$ variables over $\mathbb{C}.$ The elements in the sum on the right hand side are rational functions, which will sum together to produce an function with no poles. \newline\newline

Define local coordinates near the $i$-th fixed point $(0: ... :0:1:0:...:0)$, with 1 in the $i$-th position and 0 elsewhere, by a mapping $\mathbb{CP}^n \rightarrow \mathbb{C}^n$ for which
$$
(x_0: ...: x_i: ...: x_n) \mapsto (x_0/x_i, x_1 / x_i ..., \hat{x_i}, ... , x_n/x_i )
$$ where the hat means that we skip the entry. The local $T$ action is given by 
$$ (t_0, ..., t_n) \cdot (y_0: ... :y_n) = (\frac{t_0}{t_i}y_0, ..., \hat{y_i}, ... , \frac{t_n}{t_i}y_n)$$ thus for each $y_i$ the weight of the $T$ action is $(t_j - t_i).$ The equivariant Euler class of the tangent space $e^{T}(\nu_F,t)$ is the product of the weights of the action, so 
$$e^{T}(\nu_F,t)=\prod_{\substack{j=0, \\ j\neq i}}^n(t_j - t_i).$$

Identify $\mathbb{CP}^n = \mathbb{C}^n // U(1)$ and define the moment map $\mu$ by $$\mu: \mathbb{C}^{n+1} \rightarrow (\mathbb{R})^+ \, , \, \mu ((z_0, ..., z_n)) = (|z_0|^2, ..., |z_n|^2)$$
Suppose we have a morphisms of Lie groups $G \rightarrow H$, which gives maps $\mathbf{g} \rightarrow \mathbf{h}$ and $\mathbf{g}^* \rightarrow \mathbf{h}^*$. Let $(M, \omega)$ be a symplectic manifold upon which there is defined a Hamiltonian $G$-action with moment map $\mu_G: M \rightarrow \mathbf{h}^*$. Then there is an induced $H$ action with moment map $\mu_H : M \xrightarrow{\mu_G} \mathbf{h}^* \rightarrow \mathbf{g}^* $. For $U(1)$ the moment map is a projection 
$$\mu_{U(1)} = Pr(|z_0|^2, ..., |z_n|^2) = \sum_{i=0}^n |z_i|^2$$
Thus for the $T$-action, when evaluated at the $i-$th fixed point $f_i$ we have,
$$\mu|_{f_i}(t) = t_i.$$
Thus for each $i$ the numerator is $t_i^n$.
 The localization formula gives 
$$\int_{\mathbb{CP}^n}([\omega] +\mu(t))^n = \sum_{i=0}^n \frac{(t_i)^n}{\prod\limits_{\substack{j=0 \\ j\neq i}}^n(t_j - t_i)} \quad (*)
$$
Consider the determinant of the $(n+1)\times (n+1)$ Vandermonde matrix:
$$
\begin{vmatrix}
    1   & 1  & \dots & 1 & \dots  & 1 \\
    t_0 & t_1 & \dots & t_i & \dots  & t_n \\
    t_0^2 & t_1^2 & \dots & t_i^2 & \dots  & t_n^2 \\
    \vdots & \vdots & \vdots & \ddots & \vdots \\
    t_0^n & t_1^n & \dots & t_i^n & \dots  & t_n^n 
\end{vmatrix}
= \prod_{1\leq i < j \leq n} (t_j - t_i) \quad (* *)
$$

We may rearrange the expression in the localization formula to make the denominator Vandermonde determinant. Since in the localization formula the value of $i$ is fixed for each term in the sum, but in the determinant $i<j$ we must include a factor of $(-1)^i$ and cancel out the terms in $(**) $ that do not appear in $(*)$. Thus we may write
$$
\sum_{i=0}^n \frac{(t_i)^n}{\prod_{j=0, j\neq i}^n(t_j - t_i)}  = \frac{\sum_{i=0}^n t_i^n \cdot (-1)^i \cdot \Pi_{1\leq k < j \leq n, k\neq i}(t_k - t_j)}{\prod_{1\leq i < j \leq n} (t_j - t_i)} \,.
$$
Notice that the numerator on the right is also the determinant of an $(n+1)\times (n+1)$ Vandermonde matrix where the expansion has been done along the bottom row, so the numerator and denominator are equal. Thus we have 
$$ \int_{\mathbb{CP}^n} c^{T}_1(\mathcal{O}(1))^n =1 \,. \quad\square $$

\section{Acknowledgement}
We would like to express our deepest gratitude to our advisor, Dr. Johan Martens, for his invaluable guidance, support, and encouragement throughout this research project. His expertise, patience, and insightful feedback have been instrumental in shaping our understanding and refining our work. We are truly grateful for his mentorship and for pushing us to strive for excellence. \newline 

\noindent Furthermore, we would like to express our gratitude to the Edinburgh Mathematical Society for their generous support, which has enriched our mathematical education and facilitated the expansion of our professional network. \newline 

\newpage


\printbibliography

\end{document}